\documentclass[11pt]{amsart}
\usepackage{amsmath}
\usepackage{amsfonts}

\usepackage{amscd}
\usepackage{amsthm}
\usepackage{amssymb} \usepackage{latexsym}
\usepackage{eufrak}
\usepackage{euscript}
\usepackage{epsfig}
\usepackage{graphics}
\usepackage{array}
\usepackage{enumerate}
\usepackage{dsfont}
\usepackage{color}
\usepackage[colorlinks,citecolor=red,pagebackref,hypertexnames=false, breaklinks]{hyperref}

\newcommand{\R}{{\mathbb R}}

\newcommand{\Z}{{\mathbb Z}}

\newcommand{\Deg}{{\mathrm{Deg}}}
\newcommand{\IS}{{\mathrm{IS}_2}}

\newtheorem{lemma}{Lemma}[section]

\newtheorem{theorem}[lemma]{Theorem}
\newtheorem{corollary}[lemma]{Corollary}

\newtheorem{remark}{Remark}

\newtheorem{problem}{Problem}

\numberwithin{equation}{section}

\DeclareFixedFont{\Acknowledgment}{OT1}{cmr}{bx}{n}{14pt}
\begin{document}
\title{A note on Liouville	
 type equations on graphs}
\date{}

\author{Huabin Ge}
\address{Huabin Ge, Department of Mathematics, Beijing Jiaotong University, Beijing 100044, P.R. China}
\email{hbge@bjtu.edu.cn}


\author{Bobo Hua}
\address{Bobo Hua, School of Mathematical Sciences, LMNS, Fudan University, Shanghai 200433, China}
\email{bobohua@fudan.edu.cn}

\author{Wenfeng Jiang}
\address{Wenfeng Jiang, School of Mathematics (Zhuhai), Sun Yat-Sen University, Zhuhai, China}
\email{wen\_feng1912@outlook.com}

\maketitle

\begin{abstract}
In this note, we study the Liouville	
 equation $\Delta u=-e^u$ on a graph $G$ satisfying certain isoperimetric inequality. Following the idea of W. Ding, we prove that there exists a uniform lower bound for the energy, $\sum_G  e^u,$ of any solution $u$ to the equation. In particular, for the 2-dimensional lattice graph $\mathds{Z}^2,$ the lower bound is given by $4.$
\end{abstract}

\section{Introduction}
The Liouville equation
\begin{equation}\label{Liouvillecont}
\Delta u+e^u=0
\end{equation}
on 2-dimensional manifolds has been extensively studied in the literature. From the point of view of the theory of partial differential equations, it is critical, i.e. on the borderline of Sobolev embedding theorems in 2-dimensional case, which it is closely related to so-called Moser-Trudinger inequalities, see e.g. \cite{W}, \cite{WB}, \cite{Y} for references.

Let $u$ be a solution to the Liouville equation on the plane with finite energy, i.e.
\begin{equation}\label{eq:eq1}
\left\{
\begin{aligned}
 &  \Delta u+e^u=0\\
 &  \int_{\mathbb{R}^2}e^u<\infty.\\
 \end{aligned}
\right.
\end{equation}
An interesting argument initiated by Weiyue Ding, see \cite{CL}, shows that
$$
 \int_{\mathbb{R}^2}e^u\geq 8\pi.
$$
The key ingredient of the proof is the following isoperimetic inequality:  for any bounded domain $\Omega$ of finite perimeter in $\mathbb{R}^2,$
\begin{equation}\label{iso}
\mathrm{Length}(\partial \Omega)^2\geq 4\pi\cdot \mathrm{Area}(\Omega),
\end{equation} where $\mathrm{Length}(\partial \Omega)$ ($\mathrm{Area}(\Omega)$ resp.) denotes the length of the boundary of $\Omega$ (the area of $\Omega$ resp.).  The estimate is sharp since one can construct a family of explicit solutions,
\begin{equation}\label{eq:eq2}F_{x_0,\lambda}(x):=\ln\left[\frac{32\lambda^2}{(4+\lambda^2|x-x_0|^2)^2}\right],\quad \lambda>0,x_0\in \ \R^2,\end{equation} whose energy attain the above lower bound. Based on a delicate argument using moving plane methods, Chen and Li \cite{CL} further proved that all solutions to \eqref{eq:eq1} are exactly given by \eqref{eq:eq2}.


As is well-known, one of the difficulties for the analysis on graphs lies in the lack of chain rules for discrete Laplace operators. While linear equations have been studied extensively on graphs, people began to consider nonlinear problems on graphs such as semilinear equations recently. For semilinear equations with the nonlinearity of power type, one refers to e.g. \cite{GLY-1,GLY-2, LW-1, LW-2} .  A class of semilinear equations with the exponential nonlinearity, so-called Kazdan-Warner equations, have been studied by \cite{GLY,Ge17,GeJiang17,KS} on graphs.
The exponential nonlinearity usually causes additional difficulties for the analysis in the discrete setting. In this paper, we study the Liouville type equations on graphs, analogous to \eqref{Liouvillecont}, which are special cases of Kazdan-Warner equations. Following W. Ding's idea, we prove a uniform lower bound of the energy for the solutions to the Liouville equations on graphs satisfying isoperimetric inequalities analogous to \eqref{iso}, see Theorem~\ref{thm:main}. As a corollary, for the 2-dimensional lattice graph which is a discrete analog of $\R^2,$ we obtain an explicit lower bound for the energy of solutions to the Liouville equation, see Corollary~\ref{coro1}. This could be regarded as a preliminary step to understand the Liouville type equations on infinite graphs. 

The paper is organized as follows: In the next section, we introduce some basic setting and state our main results. Section~\ref{sec:proof} is devoted to the proof of Theorem~\ref{thm:main}.

\section{Basic setting and main results}
Let $(V, E)$ be a simple, undirected and locally finite graph, where $V$ denotes the set of vertices and $E$ denotes the set of edges.  Two vertices $x$ and $y$ are called neighbors, denoted by $x\sim y,$ if there is an edge connecting them, i.e. $\{x,y\}\in E.$ We assign weights on vertices and edges as follows:
$$
\mu:V\to (0,\infty),\quad V\ni x\mapsto \mu_x
$$
and
$$
w:E\to (0,\infty),\quad E\ni\{x,y\}\mapsto w_{xy}=w_{yx}
$$
and call the quadruple $G=(V,E,\mu,w)$ a \emph{weighted graph}. For discrete measure spaces $(V,\mu)$ and $(E,w),$ we write $\mu(A):=\sum_{x\in A} \mu_x$ and $w(B):=\sum_{e\in B} w_e$ for any subsets $A\subset V, B\subset E.$ 
For simplicity, for a function $u$ on $V$ we write
$$
\int_{V}u=\sum_{x\in V}u(x)\mu_x,
$$ whenever it makes sense.

The Laplacian on $G=(V,E,\mu,w)$ is defined as, for any function $u$ on $V$ and $x\in V,$
$$
\Delta u(x)=\frac{1}{\mu_x}\sum_{y\in V:y\sim x}w_{xy}(u(y)-u(x)).
$$ For any vertex $x,$ its weighted degree is given by $$\Deg(x):=\frac{\sum_{y:y\sim x}w_{xy}}{\mu_x}.$$ The Laplacian is bounded operator on $\ell^2(V,\mu),$ i.e. the Hilbert space of $\ell^2$ summable functions on $V$ w.r.t. the measure $\mu$, if and only if
\begin{equation}\label{bdd}\tag{$BLap$}\Deg(G):=\sup_{x\in V}\Deg(x)<\infty.\end{equation} In this paper, we always assume \eqref{bdd} holds.

For any finite subset $\Omega$ in $V,$ we denote by $$\partial \Omega:=\{\{x,y\}\in E: x\in \Omega, y\in V\setminus \Omega, \mathrm{\ or\ vice\ versa}\}$$ the (edge) boundary of $\Omega.$
We say that a weighted graph $G=(V,E,\mu,w)$ satisfies 2-dimensional isoperimetric inequality, denoted by $\IS,$ if
\begin{equation}\label{isop}\tag{$\IS$}C_{IS}:=\inf\frac{(w(\partial \Omega))^2}{\mu(\Omega)}>0,\end{equation} where the infimum is taken over all finite $\ \Omega\subset V,$ see \cite{Woe}.

In this note, we study the discrete Liouville equation
\begin{equation}\label{l-1}
  \Delta u+e^u=0.
\end{equation}
on a weighted graph $G$. Following W. Ding, see Lemma~1.1 in \cite{CL}, we obtain our main result, a discrete analog of energy estimate for the solutions to Liouville equation under the assumption of the isoperimetric inequality.
\begin{theorem}\label{thm:main}
Let $G$ be a weighted graph satisfying \eqref{bdd} and $\inf_{x\in V}\mu_x>0.$ Suppose that \eqref{isop} holds, then for any solutions $u$ of \eqref{l-1},
$$
\int_V e^u\geq \frac{C_{IS}}{\Deg(G)}.
 $$
\end{theorem}
\begin{remark}One may generalize the result to the following equation
\begin{equation*}\label{l-2}
  \Delta u+F(u)=0,
\end{equation*}
for some nonnegative function $F$ on $\mathbb{R}$ satisfying $F'\geq0$ and $F''\geq0.$
\end{remark}


We denote by $\Z^2$ the standard lattice graph with the set of vertices $\{(x,y)\in \R^2: x,y\in \Z\}$ and the set of edges $$\{\{(x_1,y_1),(x_2,y_2)\}: |x_1-x_2|+|y_1-y_2|=1\}$$ and with weights $\mu\equiv 4$ and $w\equiv1.$ It is known that it satisfies \eqref{isop} with $C_{IS}=4,$ see Theorem 6.30 in \cite{LP}. Then by the above theorem we have the following corollary.
\begin{corollary}\label{coro1}
For any solution $u$ of \eqref{l-1} on the lattice $\Z^2,$ we have
$$
\int_{\Z^2} e^u\geq 4.
 $$
\end{corollary}

This suggests the following interesting problems for further investigation.
\begin{problem} What is the sharp constant in Corollary~\ref{coro1},
i.e. $$C:=\inf \int_{\Z^2} e^u,$$ where the infimum is taken over all solutions to \eqref{l-1} on $\Z^2$?
\end{problem}

\begin{problem} Is there any solution $u$ to \eqref{l-1} on $\Z^2$ with finite energy, i.e. $\int_{\Z^2} e^u<\infty$?
\end{problem}

\section{Proof of Theorem~\ref{thm:main}}\label{sec:proof}
For any $\sigma\in\R,$ set
$$
\Omega_{\sigma}=\{x\in V|u(x)\geq \sigma\}.
$$
It is no restriction to assume that $\Omega_{\sigma}$ is finite for any $\sigma,$ otherwise by $\inf_{x\in V}\mu_x>0,$ $$\int_V e^u=+\infty.$$ By \eqref{l-1},
\begin{align*} 
\int_{\Omega_{\sigma}}e^u=&\int_{\Omega_{\sigma}}-\Delta u=\sum_{x\in\Omega_{\sigma}}\sum_{y\in V:y\sim x}w_{xy}(u(x)-u(y))\\
=&\sum_{x\in\Omega_{\sigma}}\sum_{y\in \Omega_{\sigma}: y\sim x}w_{xy}(u(x)-u(y))+\sum_{x\in\Omega_{\sigma}}\sum_{y\not\in \Omega_{\sigma}: y\sim x}w_{xy}(u(x)-u(y)) .
\end{align*}
We denote the first summand by $A.$ Then
\begin{align*}
A=&\sum_{x,y\in \Omega_{\sigma}:x\sim y}w_{xy}(u(x)-u(y))\\
=&-\sum_{x,y\in \Omega_{\sigma}:x\sim y}w_{xy}(u(y)-u(x))\\
=&-\sum_{y\in\Omega_{\sigma}}\sum_{x\in \Omega_{\sigma}: x\sim y}w_{xy}(u(y)-u(x))=-A.
\end{align*}
This yields that $A=0$ and we get
\begin{equation}\label{proof-0}
\int_{\Omega_{\sigma}}e^u=\sum_{e=\{x,y\}\in E,u(x)<\sigma\leq u(y)}w_{xy}(u(y)-u(x)).
\end{equation}
 For any $\sigma\in\R,$ let
$$
G(\sigma)=\sum_{e=\{x,y\}\in E,u(x)<\sigma\leq u(y)}w_{xy}/\big(u(y)-u(x)\big).
$$ For any subset $K\subset \R,$ we denote by $\mathds{1}_K$ the characteristic function on $K,$ i.e. $\mathds{1}_K(\sigma)=1$ if $\sigma\in K,$ and $\mathds{1}_K(\sigma)=0$ otherwise.
We have

\begin{eqnarray*}
\int_{-\infty}^{+\infty}e^{\sigma}G(\sigma)d\sigma&=&\int_{-\infty}^{+\infty}e^{\sigma}\sum_{e=\{x,y\}\in E,u(y)>u(x)}w_{xy} \big(u(y)-u(x)\big)^{-1} \mathds{1}_{(u(x),u(y)]}(\sigma)d\sigma\notag\\
&=&\sum_{e=\{x,y\}\in E,u(y)>u(x)}w_{xy}\big(u(y)-u(x)\big)^{-1}\int_{-\infty}^{+\infty}e^{\sigma} \mathds{1}_{(u(x),u(y)]}(\sigma)d\sigma \notag\\
&=&\sum_{e=\{x,y\}\in E,u(y)>u(x)}w_{xy}\frac{e^{u(y)}-e^{u(x)}}{u(y)-u(x)}\notag\\
&\leq&\sum_{e=\{x,y\}\in E,u(y)>u(x)}w_{xy}   e^{u(y)}\\
&\leq& \Deg(G) \sum_{y\in V}e^{u(y)}\mu_y,
\end{eqnarray*} where we have used the elementary inequality $\frac{e^{b}-e^{a}}{b-a}\leq e^b$ for any $a<b$ and the definition of $\Deg(G)$ in \eqref{bdd}. Hence by the above inequality,
\begin{equation}\label{square-int}
\int_{-\infty}^{+\infty}e^{\sigma}G(\sigma)\int_{\Omega_{\sigma}}e^u d\sigma\leq \int_{V}e^u\int_{-\infty}^{+\infty}e^{\sigma}G(\sigma)d \sigma\leq \Deg(G)\left(\int_{V}e^u\right)^2.
\end{equation}

On the other hand, by \eqref{proof-0} and the Cauchy-Schwarz inequality,

\begin{eqnarray*}
&&G(\sigma)\int_{\Omega_{\sigma}}e^u\\&=&\left(\sum_{e=\{x,y\}\in E,u(x)<\sigma\leq u(y)}\frac{w_{xy}}{u(y)-u(x)}\right)\left( \sum_{e=\{x,y\}\in E,u(x)<\sigma\leq u(y)}w_{xy}(u(y)-u(x))\right)\notag\\
&\geq&\big(\sum_{e=\{x,y\}\in E,u(x)<\sigma\leq u(y)}w_{xy}\big)^2=(w(\partial \Omega_\sigma))^2\notag\\
&\geq&C_{IS}\cdot\mu(\Omega_{\sigma}),\notag
\end{eqnarray*} where the last inequality follows from the isoperimetric inequality. This yields that
$$
\int_{-\infty}^{+\infty}e^{\sigma}G(\sigma)\int_{\Omega_{\sigma}}e^u\geq C_{IS}\int_{-\infty}^{+\infty} \mu(\Omega_{\sigma})e^{\sigma}=C_{IS} \int_V e^u.
$$
We prove the theorem by combining the above inequality with \eqref{square-int}.



\bigskip

{\textbf{Acknowledgements}:} We thank the anonymous referee for his/her valuable comments and suggestions. 

The research is supported by the National Natural Science Foundation of China (NSFC) under grants no. 11501027 (H. Ge) and no. 11401106 (B. Hua).





\end{document}